\documentclass[regno, 12pt]{amsart}
\usepackage{graphicx,xcolor,subcaption} \usepackage{amsmath}
\usepackage{framed,comment,enumerate}
\usepackage{amssymb,amsmath,amsfonts,amsthm,esint,enumitem,comment,mathtools,centernot}

\usepackage{hyperref}

\usepackage[style=alphabetic]{biblatex}
\renewbibmacro*{in:}{%
  \ifentrytype{article}{}{\printtext{\bibstring{in}\intitlepunct}}}

\numberwithin{equation}{section}

\newtheorem{theorem}{Theorem}[section]
\newtheorem{corollary}[theorem]{Corollary}

\theoremstyle{definition}

\DeclareMathOperator{\R}{\mathbb{R}}

\title{Non-sufficiency of smoothness in the gradient conjecture}

\author{Florian Gruen}
\address{Department of Mathematics, Graduate School of Science, Kyoto University, Kyoto, Japan}
\email{gruen.florian.32r@st.kyoto-u.ac.jp}



\begin{document}
\maketitle

\begin{abstract}
It is well known that for analytic cost functions, gradient flow trajectories have finite length and converge to a single critical point.
    The gradient conjecture of R. Thom states that, again for analytic cost functions, whenever the gradient flow trajectory converges, the limit of its unit secants exists.
    One might think that already the convergence of the gradient flow trajectory to a critical point is enough to ensure that the unit secants have a limit, but this does not hold in general - the gradient conjecture is to a certain extend sharp.\\
    We provide a counterexample in case of the missing analyticity assumption, that is a smooth (non-analytic) cost function $f$, where the limit of unit secants does not exist. In addition, $f$ satisfies even a strong geometric length-distance convergence property.
\end{abstract}

\section{Introduction}
Let $f:\R^n \rightarrow \R$ and consider the Euclidean gradient flow given by
\begin{equation}
\label{eq_continuousgradientdescent}
    x'(t) = -\nabla f (x(t)), \quad t\geq 0,\quad x(0) = x_0.
\end{equation}
In 1984, Lojasiewicz proved in a landmark paper \cite{lojagradientflow} that under the assumption $f$ being analytic, the gradient trajectory either has finite arc-length (i.e.\ it converges to a limit point $x_*$) or diverges, for $t \rightarrow \infty$. 
Recent work has been done to achieve similar results for discrete gradient decent systems and related optimization schemes, see \cite{absil}.
It is well known that merely smoothness of the cost function $f$ is not sufficient to guarantee the same result. An example are the so called "Mexican hat" functions, see \cite{absil}, \cite[page 13]{pellismelo} and  \cite{Curry1944TheMO}, that construct a function in polar coordinates, where the gradient trajectory has the whole unit circle as its limit set. 

One is not only interested whether the trajectory converges, but also in its geometry. In 1989, René Thom proposed in \cite{Thom1989ProblmesRD} the gradient conjecture. He conjectured that a convergent trajectory "does not spiral around its limit point". 
\begin{theorem}
\label{gc}
Let $x(t)$ be a solution of \eqref{eq_continuousgradientdescent} with $f$ analytic and
suppose that $x(t)\rightarrow x_*$. Then $x(t)$ has a tangent at $x_*$, that is the limit of secants $\lim_{t\rightarrow \infty}\frac{x(t)-x^*}{|x(t)-x^*|}$ exists.
\end{theorem}
Ten years later the conjecture was proven to be true in \cite{kpm}, making use of techniques from subanalytic geometry originating in \cite{lojagradientflow} and \cite{loja1}. 

We now ask the question whether smoothness could be enough to have existence of the limit of secants for a converging gradient trajectory $x$. As it turns out, the answer is no and we construct explicitly the orbit of a gradient trajectory to a smooth function $f$, where the limit of unit secants does not exists.

Upon reparametrizing by arc-length $s$, without loss of generality we assume that $x_*$ is the origin and $f(x_*)=0$, i.e.\ the differential equation in \eqref{eq_continuousgradientdescent} becomes $x'(s) = \frac{\nabla f(x(s))}{|\nabla f (x(s))|}$ and we consider $s\rightarrow s_0 < \infty$. 
In the process of the proof in \cite{kpm} a related weaker statement (than the convergence of unit secants) is also shown.
\begin{corollary}\textnormal{(``length-distance convergence") \cite[Corollary 6.5]{kpm}}
    Let $x(s)$ be a trajectory of $\frac{\nabla f}{|\nabla f|}$ and $x(s) \rightarrow 0$ as $s\rightarrow s_0$. Denote by $\sigma(s)$ its arc-length to the origin, then
    \begin{equation*}
        \frac{\sigma(s)}{|x(s)|} \rightarrow 1 \qquad \text{as } s\rightarrow s_0.
    \end{equation*}
\end{corollary}

The counterexample given here also satisfies the length-distance convergence, i.e.\ its arc-length approaches its distance to the origin. 
Morally speaking, the geometry of the trajectory resembles more and more the straight line segment connecting $x(s)$ and the origin, yet it does spiral around its limit point. For $\theta(t)$, the angle of the unit secants with $(1,0) \in \mathbb{R}^2$, not only does the limit $\lim_{t\to \infty} \theta(t)$ not exist (which would also be the case if the trajectory were to wiggle within a certain sector), but even $\lim_{t\to \infty} \theta(t)=\infty$.

\section{Counterexample}
The function is inspired by the Mexican hat function \cite{absil} and also employs the fact that $e^{-x}$ decays faster than any polynomial.

The main idea is to construct a spiral curve $\gamma: [2,\infty) \rightarrow \R^2$, written in polar coordinates as $\gamma(t)=(r(t),\theta(t))$, satisfying the desired properties and then construct a function $f(r,\theta)$ such that a gradient flow of $f$ starting at e.g.\ $x_0 = \gamma(2) = (1/2, \log(\log(2)))$  follows the trajectory of $\gamma$.

\subsection{The spiraling trajectory curve}
 Define the curve $\gamma:\R \rightarrow \R^2$ as
\begin{equation*}
    \gamma(t) = (r(t), \theta(t)) := \left(\frac{1}{t}, \log(\log(t)) \right).
\end{equation*}
We obtain the following results:
\begin{itemize}
    \item $\gamma(t)\rightarrow 0 $ as $t \rightarrow \infty$, since $r(t)=\frac{1}{t} \rightarrow 0$.
    \item $\gamma$ has finite length: 
    Since $\gamma'(t) = (x(t),y(t))' = (r(t)\cos\theta(t), r(t)\sin\theta(t))'$,
    \begin{equation*}
    \begin{split}
        \|\gamma'(t)\| &= \left( (r'\cos(\theta) - r \sin (\theta) \theta')^2 + (r' \sin (\theta) + r \cos (\theta) \theta')^2 \right)^{1/2} \\
        &=\sqrt{r'^2 \cos^2\theta -2r' r \cos\theta \sin\theta \theta' + r^2 \sin^2\theta \theta'^2+ r'^2 \sin^2\theta + 2r r' \sin\theta \cos\theta \theta' + r^2 \cos^2\theta \theta'^2 } \\
        &=\sqrt{r'(t)^2 + r(t)^2 \theta'(t)^2} = \sqrt{\frac{1}{t^4} + \frac{1}{t^4} \frac{1}{\log^2(t)}}.
    \end{split}
    \end{equation*}
    Therefore,
    \begin{equation*}
        \int_2^\infty \|\gamma'(t)\|dt  = \int_2^\infty \sqrt{\frac{1}{t^4} + \frac{1}{t^4} \frac{1}{\log^2(t)}}dt \leq C \int_2^\infty \frac{1}{t^2} < \infty.
    \end{equation*}
    \item length-distance convergence, i.e.\ $\frac{\sigma(s)}{|\gamma(s)|}=\frac{\int_s^\infty \|\gamma'(t)\|dt}{|\gamma(s)|} \rightarrow 1$ as $s\rightarrow \infty$: For this, we estimate
    \begin{equation*}
    \begin{split}
        \int_s^\infty \|\gamma'(t)\|dt  &=  \int_s^\infty \sqrt{\frac{1}{t^4} + \frac{1}{t^4} \frac{1}{\log^2(t)}}dt \\ 
        &\leq \int_s^\infty \frac{1}{t^2} + \int_s^\infty \frac{1}{t^2} \frac{1}{\log^2(t)} dt  \\
        &\leq \frac{1}{s} + \frac{1}{\log^2(s)} \int_s^\infty \frac{1}{t^2} dt \\
        &= \frac{1}{s} + \frac{1}{s}\frac{1}{\log^2(s)}. \\
    \end{split}
    \end{equation*}
    Which in turn implies
    \begin{equation*}
        \frac{\sigma(s)}{|\gamma(s)|}=\frac{\int_s^\infty \|\gamma'(t)\|dt}{|\gamma(s)|} \leq \frac{\frac{1}{s} + \frac{1}{s}\frac{1}{\log^2(s)}}{\frac{1}{s}} = 1 + \frac{1}{\log^2(s)} \rightarrow 1,
    \end{equation*}
    and thus $1 \leq \lim_{s \rightarrow \infty} \frac{\sigma(s)}{|\gamma(s)|}=1$.
    \item However, the unit secants, which are given in polar coordinates as 
    $$\frac{\gamma(t) - 0}{|\gamma(t)-0|}=(1,\theta(t))=(1,\log(\log(t))),$$ do not converge.
\end{itemize}
To conclude, $\gamma(t)$ converges to the origin, its arc-length approaches its distance to the origin, yet it spirals and hence its unit secants do not converge.

\subsection{The cost function}
It remains to construct a function $f$, having exactly $\gamma$ as a gradient flow trajectory.
 We define the function $f:B_{1/2} \rightarrow \R$ as
\begin{equation*}
\begin{split}
    f(r,\theta) &\coloneqq e^\frac{-1}{r}\left( 1- a(r) \sin(\theta - \log(\log(1/r))) \right) \\
    &=e^\frac{-1}{r}\left( 1- \frac{ \log (1/r)}{1 + r^2 \log(1/r)^2} \sin\left(\theta - \log(\log(1/r))\right) 
    \right),
\end{split}
\end{equation*}
for the auxiliary function $a(r) = \frac{\log(1/r)}{1+r^2\log(1/r)^2}$.
The exponential factor gives smoothness of the function and all its derivatives at the origin. Actually, $D^\alpha f(0,\theta)=0$ for any multi-index $\alpha$, as the factor $e^{\frac{-1}{r}}$ decays faster than any polynomial expression in $\frac{1}{r}$ and $\log(\frac{1}{r})$ increases when $r \rightarrow 0$.
By e. g. Whitney extension theorem the function could be extended to $\Tilde{f}$ defined on $\R^2$, here it is sufficient to consider the behaviour close to the origin. Note its similarity to \cite[Equation (8)]{absil}, the fractional scaling prefactor $a(r)$ and the use of $\sin$, vanishing if $(r,\theta)$ lie on $\gamma$.

We show that the gradient flow of $f$, starting at a point $x_0$ on $\gamma$, stays always on $\gamma$, i.e.\ the polar gradient $(\frac{\partial f}{\partial r}, \frac{\partial f}{ \partial \theta} \frac{1}{r})$ at points of the form 
$$\{(r,\theta): \theta - \log(\log(1/r)) = 2\pi k \quad \text{for some $k$ in $\mathbb{Z}$}\}$$
(which are precisely the points on $\gamma$) is parallel to $-\gamma'$. 
Firstly, we note
\begin{equation*}
    \gamma'(t) = (r'(t), \theta'(t)) = \left( \frac{-1}{t^2}, \frac{1}{\log(t)}\frac{1}{t} \right) =  \left(-r(t)^{2}, \frac{r(t)}{\log(1/r(t))} \right).
\end{equation*}
It remains to check that for a scaling function $b(r)$ (in fact $b(r) = \frac{1}{r^2} e^{\frac{-1}{r}} \left( \frac{\log^2(1/r)}{1+r^2\log^2(1/r)} \right)$) on the trajectory, that is whenever $\theta = \log(\log(1/r))$, we have
\begin{equation*}
    \frac{\partial f}{\partial r} = b(r) r^{
    2} \qquad \text{and} \qquad \frac{\partial f}{ \partial \theta} \frac{1}{r} = b(r) \frac{-r}{\log(1/r)}.
\end{equation*}

\noindent
We compute that,
\begin{equation*}
\begin{split}
    \frac{\partial f}{ \partial r }\bigg|_\gamma&= \frac{1}{r^2}e^\frac{-1}{r}\left( 1- a(r) \sin(\theta - \log(\log(1/r)) \right) \bigg|_\gamma  \\
    &+ e^\frac{-1}{r}\left( -a'(r) \sin(\theta - \log(\log(1/r)) -a(r) \cos(\theta - \log(\log(1/r)) \frac{1}{r^2 \log(1/r)} \right) \bigg|_\gamma\\
    &= \frac{1}{r^2} e^{\frac{-1}{r}} \left( 1- \frac{ a(r)}{ \log(1/r)}\right) \\  
    &= \frac{1}{r^2} e^{\frac{-1}{r}} \left( 1- \frac{1}{1+ r^2\log^2(1/r)}\right)\\  
    &= e^{\frac{-1}{r}} \left( \frac{\log^2(1/r)}{1+r^2\log^2(1/r)} \right), \\
    \frac{\partial f}{\partial \theta} \bigg|_\gamma &= -e^{\frac{-1}{r}} a(r) \cos(\theta - \log(\log(1/r)) \\ &= -e^{\frac{-1}{r}} a(r) \\ &= -e^{\frac{-1}{r}} \frac{\log(1/r)}{1+ r^2 \log^2(1/r)}.
\end{split}
\end{equation*}
Thus it follows that
\begin{equation*}
\begin{split}
      b(r) r^2 &= e^{\frac{-1}{r}} \left( \frac{\log^2(1/r)}{1+r^2\log^2(1/r)} \right) = \frac{\partial f}{\partial r},  \\ 
      b(r) \frac{-r}{\log(1/r)} &= 
     \frac{-1}{r} e^{\frac{-1}{r}} \left( \frac{\log^2(1/r)}{1+r^2\log^2(1/r)} \right) = \frac{\partial f}{\partial \theta} \frac{1}{r} .
\end{split}
\end{equation*}
Thus for the complete gradient of $f$ evaluated for points on $\gamma$, we obtain
\begin{equation*}
    \begin{split}
        -\nabla_{r,\theta} f |_{\gamma} =- \left( \frac{\partial f}{ \partial r}, \frac{\partial f}{ \partial \theta} \frac{1}{r} \right)\bigg|_\gamma = b(r)\left(\frac{-1}{r^2}, \frac{r}{\log(1/r)}  \right). 
    \end{split}
\end{equation*}
Since $b(r)>0$ for all $0<r\leq \frac{1}{2}$, the gradient trajectory follows $\gamma$ in the positive direction and does not stall on its way to the origin.
Therefore $ -\nabla_{r,\theta} f |_{\gamma}$ is parallel to $\gamma'$ for any point on $\gamma$, so the gradient flow stays on the curve (albeit with a different speed than its parametrization). Moreover, since $\frac{\partial f}{\partial r}$ and $\frac{\partial f}{\partial \theta}$ converge to zero as $r\rightarrow 0$ ($a(r)\rightarrow 0$ as the exponential term dominates), the origin is a critical point.

To conclude, the trajectory of the continuous gradient descent flow converges to a critical point, has convergent arc-length, even length-distance convergence, yet the secants do not have any limit.

It would be interesting to know relations to even stronger algebraic-geometric results like the \textit{(analytic) finiteness conjecture for the gradient}, see \cite{kpm}.

\begin{figure}[htbp]
\begin{minipage}[c]{0.5\linewidth}
\includegraphics[width=\linewidth]{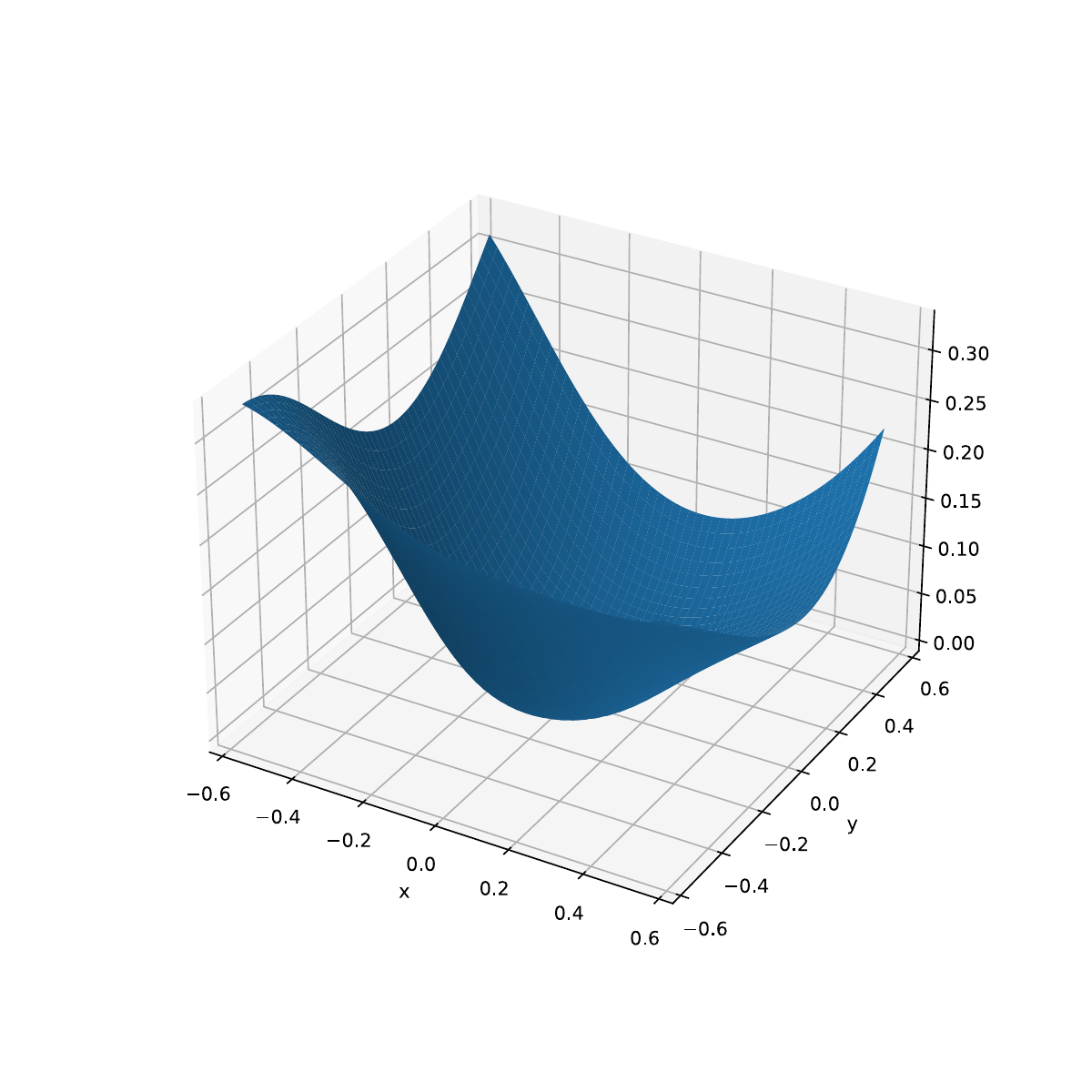}
\caption{The cost function $f$}
\end{minipage}
\hfill
\begin{minipage}[c]{0.49\linewidth}
\includegraphics[width=\linewidth]{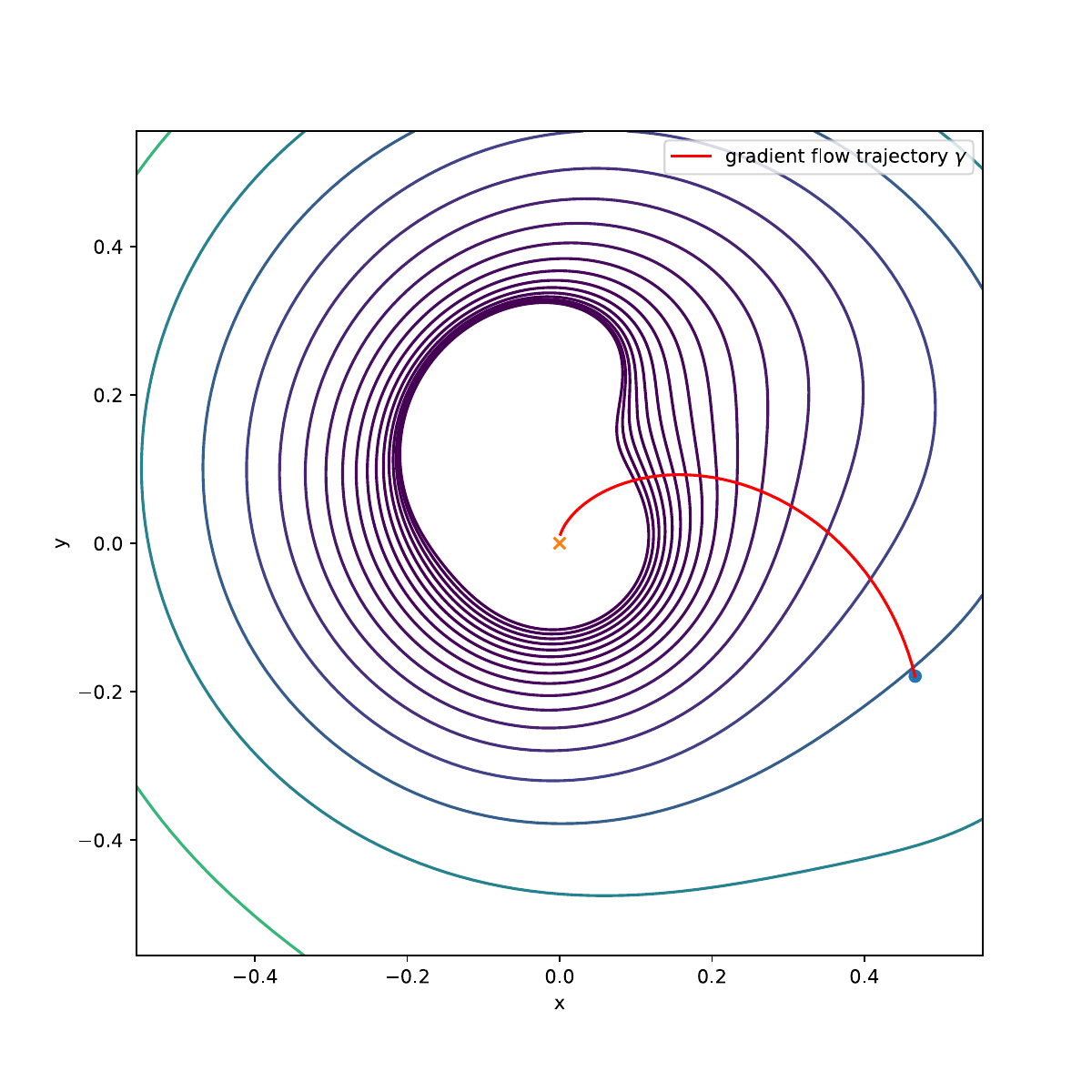}
\caption{Contour plot of $f$ and trajectory $\gamma$}
\end{minipage}%
\end{figure}

\section{Acknowledgements}
This work was produced during a semester research project at École Polytechnique Fédérale de Lausanne (EPFL) with the Chair of Continuous Optimization (OPTIM).
The author would like to thank Quentin Rebjock and Nicolas Boumal for the discussions and guidance during the project.

\printbibliography

\end{document}